\documentclass[a4paper,11pt,oneside,reqno]{amsart}
\usepackage{amsmath}
\usepackage{amsthm}
\usepackage{amssymb}
\usepackage{graphicx}
\theoremstyle{plain}
\newtheorem{teore}{Theorem}[section]
\newtheorem*{teore*}{Theorem}
\newtheorem*{coro*}{Corollary}
\newtheorem{defin}[teore]{Definition}
\newtheorem{lem}[teore]{Lemma}
\newtheorem{coro}[teore]{Corollary}
\newtheorem{propo}[teore]{Proposition}
\newtheorem{claim}{Claim}
\newtheorem*{claim*}{Claim}
\theoremstyle{remark}
\newtheorem{ejemplo}[teore]{{\sc Example}}
\newtheorem{notas}[teore]{{\sc Remark}}
\newcommand{\nrm}[1]{\|#1\|}

\newcommand{\tnrm}[1]{\|\hspace{-0.4mm}|#1\|\hspace{-0.4mm}|}

\newcommand{\prop}{\begin{propo}}
\newcommand{\fprop}{\end{propo}}
\newcommand{\cor}{\begin{coro}}
\newcommand{\fcor}{\end{coro}}

\newcommand{\defi}{\begin{defin}\rm}
\newcommand{\fdefi}{\end{defin}}
\newcommand{\eje}{\begin{ejemplo}}
\newcommand{\feje}{\end{ejemplo}}
\newcommand{\lema}{\begin{lem}}
\newcommand{\flema}{\end{lem}}
\newcommand{\teor}{\begin{teore}}
\newcommand{\fteor}{\end{teore}}
\newcommand{\nota}{\begin{notas}\rm}
\newcommand{\fnota}{ \end{notas}}
\newcommand{\clam}{\begin{claim}}
\newcommand{\fclam}{\end{claim}}
\newcommand{\clams}{\begin{claim*}}
\newcommand{\fclams}{\end{claim*}}

\newcommand{\lclam}{\begin{lclaim}}
\newcommand{\flclam}{\end{lclaim}}
\newcommand{\prucl}{\prue[Proof of Claim:]}
\newcommand{\fprucl}{\fprue}
\newcommand{\ben}{\begin{enumerate}}
\newcommand{\een}{\end{enumerate}}
\newcommand{\bit}{\begin{itemize}}
\newcommand{\eit}{\end{itemize}}

\newcommand{\mc}[1]{\mathcal{#1}}

\newcommand{\bsy}[1]{\boldsymbol{#1}}

\newcommand{\casos}{\begin{itemize}}
\newcommand{\fcasos}{\end{itemize}\setcounter{cs}{1}}

\newcommand{\ro}{\varrho}

\newcommand{\conj}[2]{ \{ {#1}\,:\,{#2} \} }

\newcommand{\ou}{\omega_{1}}

\newcommand{\om}{\omega}

\newcommand{\buit}{\emptyset}

\newcommand{\ga}{\gamma}

\newcommand{\ko}{\ensuremath{K}}
\newcommand{\co}{\ensuremath{c_{00}(\omega_1)}}

\newcommand{\al}{\alpha}
\newcommand{\be}{\beta}
\newcommand{\de}{\delta}

\newcommand{\la}{\lambda}

\newcommand{\sig}{\sigma}

\newcommand{\vep}{\varepsilon}

\newcommand{\eqs}{{\mathfrak X}}
\newcommand{\R}{{\mathbb R}}

\newcommand{\N}{{\mathbb N}}

\renewcommand{\ker}{\mathrm{Ker}}
\newcommand{\rest}{\negmedspace\negmedspace\upharpoonright\negthickspace}

\newcommand{\supp}{\mathrm{supp\, }}

\newcommand{\ran}{\mathrm{ran\, }}

\newcommand{\con}{\subseteq}

\newcommand{\prue}{\begin{proof}}
\newcommand{\fprue}{\end{proof}}

\makeindex

\begin{document}
\title{A $c_0$-saturated Banach space with no long unconditional basic sequences}
\author{J. Lopez-Abad}
\address{ Universit\'{e} Paris 7-CNRS, UMR 7056,  2 place Jussieu, 75251 Paris Cedex 05, France.}

\author{S. Todorcevic}
\address{Department of Mathematics, University of Toronto, Toronto, Canada, M5S
3G3}


\begin{abstract}
We present a Banach space $\mathfrak X$ with a Schauder basis of length $\omega_1$ which is
saturated by copies of $c_0$ and such that for every closed decomposition  of a closed subspace
$X=X_0\oplus X_1$, either $X_0$ or $X_1$ has to be separable. This can be considered as the
non-separable counterpart of the notion of hereditarily indecomposable space.  Indeed, the
subspaces of $\mathfrak X$ have ``few operators'' in the sense that every  bounded operator $T:X
\rightarrow \mathfrak{X}$ from a subspace $X$ of $\mathfrak{X}$ into $\mathfrak{X}$ is the sum of a
multiple of the inclusion and a $\omega_1$-singular operator, i.e., an operator $S$ which is not an
isomorphism on any non-separable subspace of $X$. We also show that while $\eqs$ is not distortable
(being $c_0$-saturated), it is arbitrarily $\ou$-distortable in the sense that for every $\la>1$
there is an equivalent norm $\tnrm{\cdot}$ on $\eqs$ such that for every non-separable subspace $X$
of $\eqs$ there exist $x,y\in S_X$ such that $\tnrm{x}/\tnrm{y}\ge \la$.

\end{abstract}

\maketitle
\section{Introduction}
The solutions of the unconditional basic sequence problem \cite{go-ma}
 and the distortion problem (\cite{G0}, \cite{od-schl}) have an intricate connection. They both
 have profited on one side from the development of the Tsirelson-like constructions of conditional
 norms and on the other from the development of the infinite-dimensional Ramsey theory. These
 connections are however   well understood only at the level of separable spaces. This paper,
 which can be considered as a natural continuation of our previous paper \cite{arg-lop-tod},  is an
 attempt to explain these connections in the non-separable context as well. As byproducts we
 discover some new non-separable as well as separable phenomena. Our results are all based on a
 method of constructing Banach spaces with long Schauder bases of length $\ou$, a method which among other things crucially uses the information about the
 way the classical Ramsey theorem \cite{ramsey} fails in the uncountable context \cite{tod1}.

Recall that an infinite dimensional Banach space $X$ is \emph{indecomposable} if for every closed
decomposition $X=X_0\oplus X_1$ one of the subspaces $X_0$ or $X_1$ must be finite dimensional. The
space $X$ is \emph{hereditary indecomposable}, HI in short, if every  closed infinite dimensional
subspace of $E$ is indecomposable. The first example of such a space was constructed by
Gowers-Maurey \cite{go-ma} as a byproduct  of their solution of the unconditional basic sequence
problem. The paper \cite{arg-lop-tod} considers the unconditional basic sequence problem in the
context of Banach spaces that are not necessarily separable.  In particular \cite{arg-lop-tod}
produces a reflexive Banach space $ \eqs_{\ou}$ with a Schauder basis $(e_\al)_{\al<\ou}$ of length
$\ou$ with no infinite unconditional basic sequence. Applying Gowers' dichotomy \cite{G1}, one
concludes that, while $\eqs_{\ou}$ is decomposable (as for example $\eqs_{\ou}=\overline{\langle
e_\al \rangle_{\al<\om}}\oplus \overline{\langle e_\al \rangle_{\al\ge\om}}$) it is saturated by
hereditarily indecomposable subspaces (HI-\emph{saturated} in short), i.e.,  every closed infinite
dimensional subspace of $\eqs_{\ou}$ contains an hereditarily indecomposable Banach space.

In this paper we present an example of a Banach space $\eqs$ with a Schauder basis of length $\ou$
with some extreme discrepancies  between properties of the class of separable subspaces of $\eqs$
and the class of non-separable subspaces of $\eqs$. For example, at the separable level, $\eqs$ is
\emph{saturated} by copies of $c_0$, i.e. every closed infinite dimensional  subspace of $\eqs$
contains an isomorphic copy of $c_0$. So, in particular, every closed infinite dimensional subspace
contains an infinite unconditional basic sequence. On the non-separable level, $\eqs$ contains no
unconditional basic sequence of length $\ou$. More precisely, for every closed subspace $X$ of
$\eqs$ and every decomposition $X=X_0\oplus X_1$ one of the spaces $X_0$ or $X_1$ must be
separable. In fact every bounded operator $T:X\to \eqs$, where $X$ is a closed subspace of $\eqs$,
can be decomposed as
$$T=\la i_{X,\eqs}+S,$$
where $i_{X,\eqs}:X\to \eqs$ is the inclusion mapping, and $S$ is an \emph{$\ou$-singular
operator}, the non-separable counterpart of the notion of \emph{strictly singular operator}, which
requires that $S$ is not an isomorphism on any non-separable subspace of $X$.

Another discrepancy between the behavior of separable subspaces of $\eqs$ and non-separable ones
(as well as a striking distinction between $\eqs$ and $\eqs_{\ou}$) comes when one considers the
distortion constants of its equivalent norms. Recall that the \emph{distortion} of an equivalent
norm $\tnrm{\cdot}$ of a Banach space $(X,\nrm{\cdot})$ is the constant
$$d(X,\tnrm{\cdot})=\inf_Y \sup\conj{\frac{\tnrm{x}}{\tnrm{y}}}{x,y\in S_{(Y,\nrm{\cdot})}}$$
where the infimum is taken over all \emph{infinite dimensional} subspaces $Y$ of $X$. One says that
$(X,\nrm{\cdot})$ is \emph{arbitrary distortable} if
$$\sup_{\tnrm{\cdot}}d(X,\tnrm{\cdot})=\infty$$
where the supremum is taken over all equivalent norms  $\tnrm{\cdot}$ of $(X,\nrm{\cdot})$.  The analysis of
our space $\eqs$ suggests the following variation of this notion of distortion. Given an equivalent
norm $\tnrm{\cdot}$ of a   Banach space $(X,\nrm{\cdot})$, let
$$d_{\ou}(X,\tnrm{\cdot})=\inf_Y \sup\conj{\frac{\tnrm{x}}{\tnrm{y}}}{x,y\in S_{(Y,\nrm{\cdot})}}$$
where   the infimum now is taken over all  \emph{non-separable}  subspaces $Y$ of $X$. We say that
$(X,\nrm{\cdot})$ is \emph{arbitrarily $\ou$-distortable} when
$$\sup_{\tnrm{\cdot}}d_{\ou}(X,\tnrm{\cdot})=\infty$$
where again the supremum is taken over all equivalent norms  $\tnrm{\cdot}$ of $(X,\nrm{\cdot})$. It has been
shown in \cite{arg-lop-tod} that the space $\eqs_{\ou}$ is arbitrarily distortable. Note however
that our space $\eqs$ is \emph{not} distortable at all, i.e. $d(\eqs,\tnrm{\cdot})=1$ for every
equivalent norm $\tnrm{\cdot}$ on $\eqs$. This is a consequence of the fact that $\eqs$ is
$c_0$-saturated and the well-known result of R. C. James \cite{jam} which states that if a Banach
space contains isomorphic copies of $c_0$ then it contains almost isometric copies of $c_0$.
Nevertheless, it turns out that the space $\eqs$ is distortable in the non-separable sense, i.e.,
$\eqs$ is arbitrarily $\ou$-distortable. It follows that, while the arbitrary distortion of a
Banach space $X$ implies its arbitrary $\ou$-distortion, the converse implication is not true.

\section{Definition of the space $\eqs$}

The construction of the space $\eqs$ relies on the construction of the Banach space $\eqs_{\ou}$
 from  \cite{arg-lop-tod}. So, in order to avoid unnecessary repetitions, we assume the reader is
familiar with standard definitions and results in this area (see for example \cite{lind-tza} and
\cite{arg}), and in particular with the way this has been amplified in \cite{arg-lop-tod} to the
non-separable context.

The space $\eqs$ will be defined as the completion of $(c_{00},\|\cdot\|)$ under the norm
$\|\cdot\|$ induced by a set of functionals $K  \subseteq c_{00}(\ou)$.
\defi\label{ljigrjiogj}
Let $K $ bet the minimal subset of $c_{00}(\omega_1)$ satisfying the following conditions:

\noindent {(i)} It contains $(e_{\gamma}^*)_{\gamma<\omega_1}$, is symmetric (i.e., $\phi\in K$
implies $-\phi\in K$) and is closed under the restriction on intervals of $\omega_1$.

\noindent {(ii)} For every \emph{separated} block sequence $(\phi_i)_{i=1}^{d}\subseteq \ko$, with
$d\le n_{2j}$,  one has that the  combination  $ ({1}/{m_{2j}})\sum_{i=1}^{n_{2j}}\phi_i\in\ko$.

\noindent {(iii)} For every separated \emph{special} sequence $(\phi_i)_{i=1}^{d}\subseteq \ko$
with $d\le n_{2j+1}$ one has that $\phi=({1}/{m_{2j+1}})\sum_{i=1}^{n_{2j+1}}\phi_i$ is in $K$. The
functional  $\phi$  is called a \emph{special functional}.

\noindent {(iv)} It is rationally convex.

Whenever $\phi \in K$ is of the form  $\phi=(1/m_j)\sum_{i<d}\phi_i$ given in {(ii)} or {(iii)}
above we say that $\phi$ has a \emph{weight} $w(\phi)=m_j$. Finally, the norm on $c_{00}(\omega_1)$
is defined as
$$\nrm{x}=\sup\{\phi(x)=\sum_\al\phi(\alpha)\cdot x
(\alpha):\;\phi\in\ko\}$$
 and $\eqs$ is the completion of $(\co,\nrm{\cdot})$.
\fdefi
Before we discuss the new notion of separated sequence, let us give a list of direct consequences
from the definition of $\eqs$.
\nota\label{oieutiougf}

\noindent (a)  It is clear that the norming set $K$ presented here is a subset of the one
introduced in \cite{arg-lop-tod} for the Banach space $\eqs_{\ou}=(\eqs_{\ou},\nrm{\cdot}_{\ou})$. So, it
follows that for every $x\in c_{00}(\ou)$ one has that $\nrm{x}\le \nrm{x}_{\ou}$. This fact will
be used frequently in this paper.

\noindent (b) By the minimality of $K$, there is  the following natural notion of \emph{complexity}
of every element $\phi$ of $K$: Either  $\phi=\pm e_\al^*$, or $\phi$ is a rational convex
combination $\phi=\sum_{i<k}r_i f_i$ of elements $(f_i)_{i<k}$ of $K$, or
$\phi=(1/m_j)\sum_{i<d}f_i$ for a separated block sequence $(f_i)_{i<d}$ in $K$ with $d\le n_j$.
And in this latter case we say that $w(\phi)=m_j$ is a \emph{weight} of $\phi$.

\noindent (c)  The property (i) makes  the natural Hamel basis $(e_{\al})_{\al<\omega_1}$ of
$c_{00}(\ou)$ a \emph{transfinite bimonotone Schauder basis} of $\eqs$, i.e. $(e_\al)_{\al<\ou}$ is
total and for every interval $I\con \ou$ the corresponding projection $P_I:\eqs\to
\eqs_I=\overline{\langle e_\al\rangle_{\al\in I}}$ has norm 1. Let us set $P_{\ga}=P_{[0,\ga]}$ and
$\eqs_\ga=\eqs_{[0,\ga]}$ for every countable ordinal $\ga$. It follows that every closed  infinite
dimensional  subspace contains a further subspace isomorphic to the closed linear span of a block
sequence of the basis $(e_\al)$ (see Proposition 1.3 in \cite{arg-lop-tod} for full details). This
goes in  contrast with the corresponding property of a Banach space with a Schauder basis
$(x_k)_{k\in \N}$ for which it is well-known that every closed infinite dimensional subspace
contains almost isometric copies of the closed linear span of a certain block sequence.

\noindent (d) The second property (ii) is  responsible of
 the existence of semi-normalized averages in the span of every uncountable block sequence of $\eqs $.
 The third (iii) and fourth (iv) properties  makes every operator from a closed subspace of $\eqs$ into $\eqs$ a multiple of
the identity plus a $\ou$-singular operator.

\noindent (e) The basis $(e_\al)_{\al<\ou}$  is shrinking, i.e. $(e_{\al_n})_n$ is shrinking in the
usual sense for every increasing sequence $(\al_n)_n$  of countable ordinals (the proof is
essentially equal to that for the space $\eqs_{\ou}$ provided in \cite{arg-lop-tod}; we leave the
details to the reader). It follows that $(e_\al)_{\al<\ou}$ is an uncountable weakly-null sequence,
i.e. for every $x^*\in \eqs^*$ the numerical sequence $(x^*(e_\al))_{\al<\ou}\in c_0(\ou)$. This
last property readily implies the following.

\noindent (f)  Suppose that $T:\eqs\to \eqs$ is a bounded operator. Then for every uncountable
subset $A$ of $\ou$  and every countable ordinal $\ga$ one has that $P_\ga(T(e_\al))=0$ for all but
countably many $\al\in A$.

\fnota

Fix from now on a  function $\varrho:[\ou]^2\to \omega$ with the following properties:

\noindent (i) $\varrho(\al,\ga)\le \max\{\varrho(\al,\be),\varrho(\be,\ga)\}$ for all
$\al<\be<\ga<\ou$.

\noindent (ii) $\varrho(\al,\be)\le \max\{\varrho(\al,\ga),\varrho(\be,\ga)\}$ for all
$\al<\be<\ga<\ou$.

\noindent (iii) $\conj{\al<\be}{\varrho(\al,\be)\le n}$ is finite for all $\be<\ou$ and $n\in \om$.

The reader is referred to \cite{tod1} and \cite{tod2} for full discussion of this notion and
constructions of various $\ro$-functions. We shall use it here to measure ``distances'' between
various subsets of $\ou$.

\defi
For two subsets  $s$ and $t$   of $\omega_1$,  set
$$\varrho(s,t)=\min\conj{\ro(\al,\be)}{\al\in s, \, \be\in t}.$$
Given an integer $p$ we say that $s$ and $t$ are \emph{$p$-separated} if $\ro(s,t)\ge p$. A
sequence $(s_i)_i$ is $p$-separated if it is pairwise $p$-separated, i.e. $s_i$ and $s_j$ are
$p$-separated for every $i\neq j$. The sequence $(s_i)_i$ is \emph{separated} if it is $|\bigcup_i
s_i| $-separated.

Every notion we introduced here for sets of ordinals can be naturally transferred, via their
supports, to vectors of $\eqs$.
\fdefi
The following  is a simple, but  useful statement, will give use separated subsequences of any
sufficiently long sequence of finite sets. It is the kind of result eventually used in showing that
the norm of $\eqs$ keeps a substantial conditional structure when restricted on an arbitrary
non-separable subspace of $\eqs$.

\prop
Let $(A_{i})_{i<n}$ be a block sequence of subsets of $\ou$, each of them of order-type $\om$. Then
for every block sequence  $(s_\al)_{\al \in \cup_{i<n}A_i}$  of finite sets of countable ordinals
and every integer $p$ there are $\al_i\in A_i$ ($i<n$)  such that $(s_{\al_i})_{i<n}$ is
$p$-separated.
\fprop
\prue
This is done by induction on $n$. Fix all data as in the statement for $n>1$. Then let $\al_{n-1}=\min
A_{n-1}$. Since the set
$$\overline{s_{\al_{n-1}}}^p=\conj{\be<\ou}{\text{there is some $\ga\in s_{\al_{n-1}}$ with $\be\le \ga$ and
$ \varrho(\be,\ga)\le p$}}$$ is, by property (iii) of $\varrho$, finite, one easily obtains
infinite subsets $B_i$ of $A_i$ ($i<n-1$) such that  $s_{\al}\cap \overline{s_{\al_{n-1}}}^p=\buit$
for every $\al\in B_i$ and all $i<n-1$.  By inductive hypothesis, there are $\al_i\in B_i$
($i<n-1$) such that $(s_{\al_i})_{i<n-1}$ is $p$-separated. Obviously $(s_{\al_i})_{i<n}$ is the
desired $p$-separated sequence.
\fprue

\cor\label{oeuwtu3432423fgn}
Let $n$ be an integer and let  $(s_\al^i)_{\al<\ou}$ be a block sequence of finite set of countable ordinals
for every $i<n$. Then there are $\al_0<\dots <\al_{n-1}$ such that $(s_{\al_i}^i)_{i<n}$ is a separated block
sequence.
\fcor
\prue
Let $A$ be an uncountable set such that $|s_\al^i|=p_i$ for every $i\in A$ and every $i<n$. Now for
each $i<n$, let $A_i\con A$ be of order type $\om$ and such that $A_i<A_j$  and $s_\al^i<s_{\be}^j$
if $i<j<n$  and $\al\in A_i$, $\be\in A_j$ are such that $\al<\be$. Then apply the previous
proposition to $(t_\al)_{\al\in A_0\cup \dots \cup A_{n-1}}$ and $\sum_{i<n}p_i$, where
$t_\al=s_{\al}^i$ for the unique $i<n$ such that $\al\in A_i$.
\fprue

\subsection{Rapidly increasing sequences and deciding pairs}

We now  introduce some standard technical tools in this field. Particularly,  a sort of vectors,
called $\ell_1^k$-averages, and the so-called rapidly increasing sequences (RIS in short). The
importance of rapidly increasing sequences $(x_n)$ is that it is possible to estimate the norm of
linear combinations of $(x_n)$ in terms of norms of linear combinations of the basis
$(e_\al)_{\al<\ou}$. In this sense RIS behave like a subsequence $(e_{\al_k})_k$ of basis
$(e_\al)_{\al<\ou}$.  The role of  $\ell_1^k$-averages is that they are useful in creating RIS. We
give the precise definitions now.

\defi
Let $C,\vep>0$. A normalized block sequence  $(x_k)_k$ of $X$ is called  a $(C,\vep)$-
\emph{rapidly increasing sequence} ($(C,\vep)$-RIS in short) iff there is an increasing sequence
$(j_k)_k$ of integers such that for all $k$,

\noindent (i) $\nrm{x_k}\le C$,

\noindent (ii) $|\supp x_k| \le m_{j_{k+1}}\vep$ and

\noindent (iii)  For every functional $\phi\in K$ of type I  with   $w(\phi)<m_{j_k}$ one has that
$|\phi(x_k)|\le C/w(\phi)$.

Let $C>0$ and $ k\in \N$. A normalized vector $y$ is called a $C-\ell_1^k$-\emph{\emph{average}}
iff there is a finite block sequence $(x_0,\dots,x_{k-1})$ such that
  $y=(x_0+\dots + x_{k-1})/k$ and  $\nrm{x_i}\le C$.
\fdefi

First observe that  any $\omega$-subsequence $(e_{\al_n})_n$ of the basis is a  $(1,\vep)$-RIS for
every $\vep$. Note also that  it follows easily from the definition that if $(x_n)$ is a
$(C,\vep)$-RIS, then for every $\vep'>0$ there is a subsequence $(x_n)_{n\in A}$ of $(x_n)$ which
is a $(C,\vep')$-RIS.

As the norming set $K$ we are using here is not saturated by ``free'' combinations of the form
$(1/m_{2j})\sum_{i<n_{2j}}f_i$, one cannot expect that there are $\ell_1^k$-averages in the span of
arbitrary block sequence. Indeed, as Theorem \ref{ieoithie} shows, this is not the case for most of
the block sequences, since clearly if $(x_n)$ is $C$-equivalent to the $c_0$-basis, then
$\nrm{(x_0+\dots + x_k)/k}\le C/k$ for every $k$. However, the next Proposition guarantees their
existence. Its proof is the natural modification of the standard proof for the separable case,
which can be found, for example, in \cite{A-M} or in \cite{go-ma}.

\prop
\label{dfmlejtilejiew} For every   $k\in \N$ there is  $l=l(k)\in \N$ such that the following
holds: Suppose that $(x_i)_{i<l}$ is a normalized block sequence with the property that there
exists a separated sequence $(\phi_i)_{i<l}$ which is biorthogonal to $(x_i)_{i<l}$. Then $\langle
x_i \rangle_{i<l}$ contains 2-$\ell_1^k$-averages. \qed
\fprop

We recall the following is also a well known fact about $\ell_1^k$-averages, connecting them to
RIS.
\prop\label{risav}
Suppose that $y$ is a  $C-\ell_1^k$-average  and suppose that $E_0<\dots < E_{l-1}$ are intervals
with $n< k$. Then  $ \sum_{i=0}^{l-1} \nrm{E_i y}\le C(1+{2l}/k)$. As a consequence, if $y$ is a
$C-\ell_1^{n_j}$-average and    $\phi\in K$ is  with $w(\phi)<m_j$ , then $|\phi(y)|\le
{3C}/{2w(\phi)}$.

In particular, for $2-\ell_1^{n_j}$-averages we get that $|\phi(y)|\le 3/w(\phi)$ if $w(\phi)<m_j$.
\qed
\fprop

\nota\label{gheuhrufhsd}
It follows that if $\bsy{x}=(x_k)_k$ is a normalized block sequence such that each $x_k$ is a
$2-\ell_1^{n_{j_k}}$-average with $|\supp x_k|\le m_{j_{k+1}}$, then $\bsy{x}$ is a $(3,1)$-RIS. It
readily follows that  if  $(x_k)_k$ is a normalized block sequence having a separated biorthogonal
pair, then the linear span of $(x_k)_k$ contains $(3,\vep)$-RIS for every $\vep>0$.
\fnota

 Let
$<_\mathrm{antilex}$ denote the anti-lexicographical ordering  on $\N\times \ou$. Whenever we say
that a sequence $(x_k^\al)_{(k,\al)\in A}$ indexed on a subset $A$ of $\N\times \ou$ is a block
sequence we mean that $x_k^\al$ is finitely supported and that $x_{k_0}^{\al_0}<x_{k_1}^{\al_1}$
whenever $(k_0,\al_0)<_\mathrm{antilex} (k_1,\al_1)$.  In the following definitions and  lemma, we
introduce two genuinely non-separable tools. They are necessary for us because it not true in
general that a  non-separable subspace of $\eqs$  contains an almost isometric  copy of the closed
linear span of some uncountable block sequence. Recall that if $X$ is a separable Banach space with
a Schauder basis of length $\om$, the corresponding result is  true and very frequently used.

\defi Let $C,\vep>0$.
We call a  sequence $\boldsymbol{x} =(x_{k}^\al)_{(k,\al)\in \N\times\ou}$ a  \emph{long rapidly
increasing sequence} (\emph{LRIS} in short)  iff

\noindent (i) $x_{k_0}^{\al_0}<x_{k_1}^{\al_1}$ for every
$(k_0,\al_0)<_{\text{antilex}}(k_1,\al_1)$, and

\noindent (ii) The cardinality of $\supp x_k^\al$ only depends on  $k$.

\noindent (iii) There is a sequence of integers $(j_k)$ such that  $x_{k}^\al$ is a
$2-\ell_1^{n_{j_k}}$-average and $|\supp x_k^{\al}|<m_{j_{k+1}}$ for every $(k,\al)$.

\fdefi
\nota
The chosen name is because  it follows from the definition and  Remark \ref{gheuhrufhsd} that if
$f=(f_0,f_1):\N\to \N\times \om_1$ is such that
$$\text{if $k<l$  then $f_0(k)<f_0(l)$ and $f_1(k)\le f_1(l)$}$$ then
one has that $(x_{f_0(k)}^{f_1(k)})_{k\in \N}$ is a $(3,1)$-RIS.
\fnota

\defi
Given $\vep>0$ and two vectors $x,y\in c_0(\ou)$, with $x\neq 0$, we write $x<_\vep y$ to denote
that
$$\nrm{P_{\sup\supp x}(y)}<\vep.$$
By technical reasons, we   declare $0<_\vep y$ for every $y$.

We recall that for a vector $x\in c_{0}(\ou)$, one sets $\ran x=[\min \supp x,\sup\supp x]$.
\fdefi

\defi\label{decpair}
Given a bounded operator $T:X\to \eqs$, where $X$ is a closed non-separable subspace of $\eqs$, we
say that the couple $(\bsy x,\bsy y)$ is a \emph{deciding pair} for $T$ if the following holds:

\noindent (i)  $\boldsymbol{x}=(x_k^\al)_{(k,\al)\in \N\times \ou}\con X$ and
$\boldsymbol{y}=(y_k^\al)_{(k,\al)\in \N\times \ou}$   is a LRIS and $\bsy x\con X$.

\noindent (ii) $x_{k}^\al\in X$ and  $\ran y_k^\al\con \ran x_k^\al$ for every pair $(k,\al)$.

\noindent (iii)  $\sum_{k\in \N}\nrm{x_k^{\al}-y_k^\al}\le 1$ for every $\al$.

\noindent (iv) $ x_{k_0}^{\al_0}<_{\vep_{k_1}} x_{k_1}^{\al_1}$ and    $
T(x_{k_0}^{\al_0})<_{\vep_{k_1}} T(x_{k_1}^{\al_1})$ for every
$(k_0,\al_0)<_\text{antilex}(k_1,\al_1)$.

A \emph{transversal subsequence} of a double-indexed sequence $(x_{k}^{\al})_{(k,\al)}$ is a finite
subsequence of the form $(x_{k_i}^{\al_i})_{i}$ where $k_i<k_{i+1}$ and $\al_i<\al_{i+1}$ for every
$i$.
\fdefi
In other words a deciding pair is nothing else but an uncountable ordered sequence of rapidly
increasing sequences $(y_k^\al)_{k\in \N}$ which are asymptotically closed to a sequence
$(x_k^\al)_{k\in \N}$ in $X$ for which the corresponding sequence of images $(Tx_k^\al)_{k\in \N}$
is ``almost'' block ordered.

Before we prove that deciding pairs always exist, we give some explanation of this notion.

\prop\label{oehthorhgoiwehiohtg}
Suppose that $(\bsy x,\bsy y)$ is a deciding pair for $T$. Then for every $\vep>0$ and every
integer $l$ there is a transversal subsequence $\bsy z=(y_{k_i}^{\al_i})_{i<l}$ of $\bsy y$ such
that

\noindent (i) $\bsy z$ is a $(3,\vep)$-RIS.

\noindent (ii) $\bsy z$ has a biorthogonal separated block sequence in the norming set $K$.

\noindent (iii) $\sum_{i<l}\nrm{x_{k_i}^{\al_i}-y_{k_i}^{\al_i}}\le \vep$.

\fprop

\prue Fix all data.
Let $M\con \N$ be such that for every $\al<\ou$ one has that
 $(y_{k}^\al)_{k\in M}$ is a
$(3,\vep)$-RIS  and such that $\sum_{k\in M}\nrm{x_{k}^\al-y_k^\al}\le \vep$. Fix also for each
pair $(k,\al)$ a functional  $\phi_{k}^\al\in K$ such that $\phi_k^\al(y_k^\al)=1$ and with $\ran
\phi_k^\al\con \ran y_k^\al$. Now apply Corollary \ref{oeuwtu3432423fgn} to $(\supp
\phi_k^{\al})_{\al<\ou}$ ($k\in M$) to find a transversal subsequence $(\phi_{k_i}^{\al_i})_{i<l}$.
Then $(y_{k_i}^{\al_i})_{i<l}$ is the desired sequence.
\fprue

\lema\label{onrothoirgjfig}
Every bounded operator $T:X\to\eqs$ with $X$ non-separable has a deciding pair.
\flema
\prue
First we make the following approximation to the final result.
\clam
For every integer $k$ and every $\vep>0$ there are two normalized sequences $(x_\al)_{\al<\ou}$ and
$(y_\al)_{\al<\ou}$ such that

\noindent (a) $x_\al\in X$ ($\al<\ou$).

\noindent (b) $\ran y_\al\con \ran x_\al$  and  $\nrm{x_\al-y_\al}\le \vep$ ($\al<\ou$).

\noindent (c) $y_\al$ is a $2-\ell_1^k$-average and $|\supp y_\al|$ is independent of $\al<\ou$.

\noindent (d) $x_\al<_\vep x_\be$  and $Tx_\al<_\vep Ty_\al$ for every $\al<\be$.
\fclam
Let us show the desired result from this claim: Find recursively for each $k\in \N$ two  sequences
$(z_\al^k)_{\al<\ou}$ and $(t_\al^k)_{\al<\ou}$ as the result of the application of the previous
claim to the integer $n_{j_k}$ and    $\vep_k$ and where $j_k$ is chosen such that   $|\supp
t_{\al}^{k-1}|<m_{j_k}$. Finally, it is not difficult to see that  one can extract for every $k$ a
subsequence $(x_k^\al,y_k^\al)_{\al<\ou}$  of  $(z_\al^k,t_\al^k)_{\al<\ou} $ with the property
that $x_{k_0}^{\al_0}<_{\vep_{k_1}} x_{k_1}^{\al_1}$ and $T(x_{k_0}^{\al_0})<_{\vep_{k_1}}
T(x_{k_1}^{\al_1})$ for every $(k_0,\al_0)<_\mathrm{antilex} (k_1,\al_1)$.

Let us give now a proof of the claim. Fix $k$ and $\vep>0$, and  let $l$ be any integer given by
Proposition \ref{dfmlejtilejiew} when applied to $k$.   Set $\de=\vep/2l$. Now use that the bounded
operator  $U=U(T,\ga):X\to X_\ga\oplus X_{\ga}$ defined by $U(x)=(P_\ga(x),P_\ga(T(x)))$ is
$\ou$-singular (because it has separable range) to find two normalized sequences
$\bsy{z}=(z_\al)_{\al<\ou}$ and $\bsy{t}=(t_\al)_{\al<\ou}$, and a block sequence
$(F_\al)_{\al<\ou}$ of finite sets  $F_\al\con \ou$ of size $l$ such that

\noindent (e) $\bsy{z}\con X$, and $\bsy t$ is a block sequence.

\noindent (f) $\ran t_\al\con \ran z_\al$ and $\nrm{z_\al-t_\al}\le \de$  for every $\al$.

\noindent (g) $z_\al<_\de z_\be$ and $T(z_\al)<_\de T(z_\be)$ for every $\al<\be$.

\noindent (h) For every countable ordinal $\al$ there is a separated block sequence
$(f_\xi)_{\xi\in F_\al}$ that is biorthogonal to $(t_\xi)_{\xi\in F_\al}$

We observe that (h) can be achieved by a simple application of Corollary \ref{oeuwtu3432423fgn}.
 By Proposition \ref{dfmlejtilejiew}, we
can find a $2-\ell_1^k$-average $y_\al\in \langle t_\be \rangle_{\be\in F_\al}$ for each $\al<\ou$.
It is easy to see that if $x_\al$ is an arbitrary normalized vector in $\langle z_\be
\rangle_{\be\in F_\al}$ such that $\ran y_\al\con \ran x_\al$ and $\nrm{x_\al-y_\al}\le \vep$ (and
there is such vector), then the corresponding sequences $(x_\al)_{\al<\ou}$ and $(y_\al)_{\al<\ou}$
fulfill all the properties  \emph{(a)} to \emph{(d)}.
\fprue
\nota\label{ojgtiorufjnw}
It is easy to find for given two bounded operators $T_0,T_1:X\to \eqs$ a deciding pair $(\bsy
x,\bsy y)$ for, simultaneously, $T_0$ and $T_1$. This can be done simply by replacing $U(T,\ga)$
above by the mapping $\bar{U} :X\to \eqs_{\ga}^3$ defined by $x\mapsto
\bar{U}(x)=(P_\ga(x),P_\ga(T_0(x)),P_\ga(T_1(x)))$.
\fnota

\section{Main properties of the space $\eqs$}
\subsection{$c_0$-saturation}
We are ready to prove  that $\eqs$ is $c_0$-saturated.  We start with the following more
informative result.

\lema\label{orjuitorjhgfj}
Suppose that $(x_k)$ is a normalized block sequence such that
$$\lim_{k\to \infty} \nrm{x_k}_\infty=0.$$
Then $(x_k)$ has a subsequence which is 5-equivalent to the natural basis of $c_0$.
\flema
\prue
Let $(y_k)_k$  be a subsequence of $(x_k)_k$  with $\nrm{y_k}_\infty\le 1/({2^{k+2}k})$.
\clam
There is an infinite set $M$ of integers such that for every triple $(k_0,k_1,k_2)$ in $M$ one has
that
\begin{equation}
\label{wejgiiujibhjbofdhss} \max\conj{\varrho(\al,\be)}{\al\in \supp y_{k_0},\, \be\in \supp
y_{k_1}}< k_2.
\end{equation}
\fclam
\prucl We color every triple $(k_0,k_1,k_2)$ of integers by 1 if \eqref{wejgiiujibhjbofdhss} holds
and 0 otherwise. By the classical Ramsey theorem  we can find an infinite set $M$ all whose triples
are equally colored. If this color is 1, then we are done. Otherwise, suppose that this color is 0,
and let us yield a contradiction.  Fix two integers $k_0<k_1$ in $M$. Then for every $k\in M$ with
$k
>k_1$ let $\al_k\in \supp y_{k_0}$ and $\be_k\in \supp y_{k_1}$ with $\varrho(\al_k,\be_k)\ge k$.
Find an infinite $P\con N$ with $(\al_k,\be_k)=(\al,\be)$  for every $k\in P$. Then
$\varrho(\al,\be)\ge k$ for every $k\in P$, a contradiction. \fprucl
 Fix such $M$ from the claim.  We show that $(y_n)_{n\in M}$ is 5-equivalent to the natural basis
 of $c_0$. Observe that since the basis $(e_\al)_{\al}$ is bimonotone one has that $\nrm{\sum_{n\in M} a_n y_n}\ge
 \nrm{(a_n)_{n\in M}}_{\infty}$ for every sequence $(a_n)_{n\in M}$ of scalars. So it remains to show
 that $\nrm{\sum_{n\in M} a_n y_n}\le 5
 \nrm{(a_n)_{n\in M}}_{\infty}$. This is done in the next.
\clam
For every $\phi\in K$ and every sequence $(a_n)_{n\in M}$ of scalar one has that
\begin{equation}
\label{njonhnjdgfd} |\phi(\sum_n a_n y_n)|\le 5 \nrm{(a_n)_{n\in M}}_{\infty}.
\end{equation}
\fclam
\prucl Fix all data, and set $y=\sum_{n}a_n y_n $. The proof of \eqref{njonhnjdgfd} is done by
induction on the complexity of $\phi$. If $\phi=\pm e_\al^*$, the result is trivial. Suppose that
$\phi$ is rational convex combination $\phi=\sum_{i<k} c_i f_i$ of elements $(f_i)_{i<n}$ of $K$.
Then, by applying the inductive hypothesis to $f_i$'s, one has that  $|\phi(y)\le \sum_{i<k} c_i
|f_i(y)|\le 5\nrm{(a_n)}_{\infty}$. Suppose now that $\phi=(1/m_j)\sum_{i<d}f_i$, where $d\le n_j$
and  $(f_i)_{i<d}$ is a separated block sequence in $K$. Let
\begin{align*}
a=&\conj{i<d}{\supp y \cap \supp f_i\neq \buit },\text{ and for $i\in a$ let}\\
k(i)=&\min \conj{k}{\supp y_k\cap \supp f_i\neq \buit},\text{ and} \\
L=& \conj{k(i)}{i\in a}.
\end{align*}
\clam
\noindent (a) If $k<\min L$, then $\supp y_k\cap \supp f=\buit$.

\noindent (b) There is some  is some $\bar{i}\in a$   such that for every $k>\max L$, then $\supp
y_k\cap \supp f=\supp y_k\cap \supp f_{k(i)}$.

\noindent (c) For every   two consecutive $k_0<k_1$ in $K$ there is some $\bar i(k_0,k_1)\in a$
such that for every $k$ with $k_0<k<k_1$ one has that $\supp y_k \cap \supp f=\supp y_k\cap \supp
f_{\bar i(k_0,k_1)}$.
\fclam
\prucl The statement (a) is clear. It is not difficult to show that  the statements  (b) and (c)
follow from the following, also non difficult, fact: Suppose that $I$ is an interval of integers,
and suppose that there are at least two integers $i_0$ and $i_1$ in $a$ such that $\supp
f_{i_\vep}\cap \bigcup_{k\in I}\supp y_k\neq \buit$ for $\vep=0,1$. Then $I\cap L\neq \buit$.
\fprucl Now we consider two cases:

\noindent \emph{Case 1.} The cardinality of $L$ is at most two. If $L=\buit$, then $f(x)=0$.
Suppose that $L=\{\bar{k}\}$. Then  using \emph{(a)}, \emph{(b)} and the inductive hypothesis,
above one obtains that
$$|f(x)|\le \frac{1}{m_j}|f_{\bar{i}}(\sum_{k> \bar k}a_k y_k)|+|a_{\bar k}|\le 2 \nrm{(a_k)_k}_{\infty}.$$
Finally suppose that $L=\{k_0,k_1\}$ with $k_0<k_1$. Then  one has that
\begin{align*}
|f(x)|\le & |a_{k_0}|+|a_{k_1}|+\frac{1}{m_j}\left|f_{\bar i(k_0,k_1)}(\sum_{k=k_0+1}^{k_1-1}a_k
y_k) +f_{\bar i}(\sum_{k>k_1}a_k y_k) \right| \le 5 \nrm{(a_k)_k}_\infty.
\end{align*}

\noindent \emph{Case 2.} The set $L$ has cardinality at least three. Let $k_0<k_1<k_2$ be the least
three elements of $L$. Find $i_0<i_1$ in $a$ such that $k(i_\vep)=k_\vep$ for $\vep=0,1$, and then
$\al_\vep\in \supp f_{i_\vep}\cap \supp y_{k_\vep} $ for $\vep=0,1$. It follows that
\begin{equation*}
\nrm{f}_{\ell_1}\le |\supp f|\le \varrho(\al_0,\al_1)< k_2.
\end{equation*}
Hence for every $k\ge k_2$ one has that
\begin{equation}
\label{ewfmkjirjjfd}  |f(y_k)| \le \nrm{f}_{\ell_1}\nrm{y_k}_\infty \le \frac{1}{2^{k+2}}.
\end{equation}
Using the inequality \eqref{ewfmkjirjjfd}, conditions \emph{(a)}-\emph{(c)} above    and the
inductive hypothesis applied to $f_{i(i_0,i_1)}$ and $f_{i_1}$ one obtains that
\begin{align*}
|f(x)|\le & |a_{k_0}|+|a_{k_1}|+ \frac{1}{m_j}\left|f_{\bar i(k_0,k_1)}(\sum_{k=k_0+1}^{k_1- 1}a_k
y_k)
+f_{\bar i (k_1,k_2)}(\sum_{k=k_1+1}^ {k_2-1}a_k y_k)\right|+   |f(\sum_{k\ge k_2} a_ky_k)| \le \\
\le & 2 \nrm{(a_k)_k}_{\infty}+\frac{5}{2}\nrm{(a_k)_k}_\infty+\frac12\nrm{(a_k)}_\infty = 5
\nrm{(a_n)}_\infty,
\end{align*}
as desired.
\fprucl
\fprue

\teor\label{ieoithie}
$\eqs$ is $c_0$-saturated.
\fteor
\prue
Fix an closed infinite dimensional  subspace $X$ of $\eqs$. We may assume, by Remark
\ref{oieutiougf} (c), that indeed $X$ is the closed linear span of a normalized block sequence
$(x_k)$. Suppose that $c_0$ does not embed isomorphically into $X$, and let us yield a
contradiction. Consider  the norm-one operator $  \sum_k a_k x_k\in X \mapsto (a_k)_k\in c_{0}$.
This is, by hypothesis, strictly singular. So we can find a normalized block subsequence $(y_k)_k$
of $(x_k)_k$ with $\lim_{k\to \infty}\nrm{y_k}_\infty=0$. Then by Lemma \ref{orjuitorjhgfj}
$(y_k)_k$ has a subsequence equivalent to the $c_0$-basis, a contradiction.
\fprue

\subsection{Distortion}
From the previous Theorem \ref{ieoithie} one immediately obtains the following.

\cor $\eqs$ is not
distortable. \qed
\fcor
In contrast to this, we have a strong distortion phenomenon in the level of non-separable subspaces
of $\eqs$:

\teor
$\eqs$ is   arbitrarily $\ou$-distortable.
\fteor
\prue
We follow some of the ideas   used to show that $\eqs_{\ou}$ is arbitrarily distortable (see
Corollary 5.36 in \cite{arg-lop-tod}).   For $j\in \N$, and $x\in \eqs_{\ou}$, let
$$\nrm{x}_{j}=\sup \conj{\phi(x)}{w(\phi)=m_{2j}}.$$
Notice that, obviously, for every $\phi\in K$ one has that $(\phi)$ is a dependent sequence, hence
$(1/m_{2j})\phi\in K$. It follows that $\nrm{\cdot}_j\le \nrm{\cdot}\le m_{2j}\nrm{\cdot}_j$. Fix a
closed non-separable subspace $X$ of $ \eqs$.      Let $(\bsy x,\bsy y)$ be a deciding pair for the
inclusion mapping $i_{X,\eqs}:X\to \eqs$.

Now fix an integer $l$ and $\vep>0$. Use Proposition \ref{oehthorhgoiwehiohtg} to find  a transversal
subsequence  $\bsy z=(y_{k_i}^{\al_i})_{i<n_{2l}}$ of $\bsy y$  such that

\noindent (a) $\sum_{i<l}\nrm{x_{k_i}^{\al_i}-y_{k_i}^{\al_i}}\le \vep$, and

\noindent (b) $\bsy z$ is a $(3,n_{2l}^{-2})$-RIS with a biorthogonal separated block sequence
$(f_i)_{i<n_{2l}}$ in $K$. Set
$$z_l=\frac{m_{2l}}{n_{2l}}\sum_{i<l}y_{k_i}^{\al_i}, \, \phi_l=\frac1{m_{2l}} \sum_{i<l}f_{i}.$$
Then $\phi_l\in K$, and the pair $(z_l,\phi_l)$ is what we called in \cite{arg-lop-tod} (Definition
3.1, see also the proof of Proposition 4.12) a $(6,2l)$-\emph{exact pair}. So, it follows that if
$l=j$ then one has
$$1\le \nrm{z_j}_{j}\le \nrm{z_j}\le 6,$$
while if $l>j$ then
$$\text{$1\le \nrm{z_l}\le 6$ and $\nrm{z_l}_{j}\le \frac{12}{m_{2j}}$.}$$
Hence for every $l>j$ one has   the   discrepancy
\begin{equation}\label{njhdfdj}
\frac{ \nrm{ z_j/\nrm{z_j}}_{j}}{ \nrm{ z_l/\nrm{z_l}}_{j}}\ge
\frac{1/6}{12/m_{2j+1}}=\frac{m_{2j}}{72}.
\end{equation}
Since the vectors $z_k$ and $z_l$ can be found to be arbitrary close to $X$, it follows that one
can obtain a similar  inequality  to \eqref{njhdfdj} for vectors in $X$. Hence $(\eqs,\nrm{\cdot})$
is arbitrarily $\ou$-distortable.
\fprue

\subsection{Operators}
It turns out that while the space $\eqs$ is $c_0$-saturated it has the following form of
indecomposability at the non-separable level.

\defi
A  Banach space $X$ is \emph{$\ou$-indecomposable}  if  for every decomposition $X= Y\oplus Z$  one
has that either $Y$ or $Z$ is separable. We say that $X$ is \emph{$\ou$-hereditarily
indecomposable} ($\ou$-HI in short) if every subspace of $X$ is $\ou$-indecomposable.
\fdefi
This kind of indecomposability corresponds to the following notion of singularity for operators.

\defi
An operator $T:X\to Y$ is \emph{$\ou$-singular} if $T$ is not an isomorphism on any non-separable
subspace of $X$.
\fdefi

 Observe that strictly singular and separable range operators are     $\ou$-singular. While the
 strictly singular operators and operators with separable ranges form closed ideals of the Banach
algebra  $\mc L(X)$ of all bounded operators from $X$ into $X$, we do not know if, in general, this
is also the case for the family $\mc S_{\ou}(X)$ of $\ou$-singular operators of $X$. Indeed we do
not even know if $\mc S_{\ou}(X)$ is closed under sums. We shall show however that for subspaces
$X$ of our space $\eqs$ one does have the property that  $\mc S_{\ou}(X)$  form a closed ideal in
the algebra $\mc L(X)$.

Using this  notion one can have the following sufficient condition for being $\ou$-HI.

\prop\label{jyhruthr}
Suppose that $X$ has the property that for every subspace $Y$ of $X$ every bounded operator $T:Y\to
X$ is of the form $T=\la i_{Y,X}+S$ where $S$ is $\ou$-singular and $\la\in \R$. Then $X$ is
$\ou$-HI.
\fprop
\prue
Otherwise, fix two nonseparable subspaces $Y$ and $Z$ of $X$ such that $d(S_X,S_Y)>0$.  It follows
that the two natural projections $P_{Y}:Y\oplus Z\to Y$ and $P_Z:Y\oplus Z\to Z$  are both bounded.
Fix $\la\in \R$ such that $T=i_{Y,X}\circ P_Y=\la i_{Y\oplus Z,X} + S$ with $S$ $\ou$-singular.
Since $T^2=T$, we have that
\begin{equation}\label{reuturuhggh}
(\la^2-\la)i_{Y\oplus Z,X} =  ((1-2\la) i_{Y\oplus Z,X} -S)\circ S.
\end{equation}
Since it is clear that $U\circ S$ is $\ou$-singular if $S$ is $\ou$-singular, it follows from
\eqref{reuturuhggh} that $\la^2=\la$. Without loss of generality, we may assume that $\la=1$ (if
$\la=0$ we replace $Y$ by $Z$ in the preceding argument). Since $P_Y\rest Z=0$, we obtain that $S=-
i_{Z,X}$, a contradiction.
\fprue
\nota
Recall that   V. Ferenczi  has shown in \cite{fere} that if a complex Banach spaces $X$ is HI then every
operator from a subspace $Y$ of $X$ into $X$ is a multiple of the inclusion plus a strictly singular
operator. We do not  know if the analogous result is true for $\ou$-singular operators or, in other words, if
the converse implication of Proposition \ref{jyhruthr} is true in the case of complex Banach spaces.
\fnota
The main purpose of this subsection is the study  of the operator space $\mc L(X,\eqs)$, where $X$
is an arbitrary closed infinite dimensional subspace of $\eqs$.  For the next few lemmas we fix  a
bounded operator $T:X\to \eqs$ from a closed infinite dimensional subspace $X$ of $\eqs$ into the
space $\eqs$. We also fix a deciding pair $(\boldsymbol{x},\boldsymbol{y})$ for $T$ (see Definition
\ref{decpair}).
\lema\label{ojrtreuret}
For all but countably many $\al<\ou$ one has that $\lim_{k\to \infty}d(Tx_k^\al,\R x_k^\al)=0$.
\flema
\prue
Otherwise, using the property (iii) of the deciding pair $(\bsy x,\bsy y)$ and going to subsequences if
necessary, we may assume that there is $\vep>0$ such that
$$\inf_{k\in \N}d(Tx_k^\al,\R
y_k^\al)>\vep$$
 for every
countable ordinal $\al $.  Now using Hahn-Banach theorem and the fact that the norming set $K$ is
closed under rational convex combinations and restrictions on intervals we can find for every pair
$(k,\al)\in\N \times \ou$ a functional $f_{k}^\al\in K$ such that  one has that

\noindent (a)   $(f_k^{\al})_{(k,\al)\in \N\times \ou}$ is a block sequence and $\ran f_k^\al\con
\ran x_k^\al$ for every $(k,\al)$, and

\noindent (b) $| f_k^\al(y_k^\al)|\le \vep_k$ while $f_k^\al(T(x_k^\al))\ge \vep $ for every
$(k,\al)$.

Fix $j$ with $\vep m_{2j+1} >   2\nrm{T}$. Now use   Proposition \ref{oehthorhgoiwehiohtg}  to find a
sequence $(F_i^\al)_{(i,\al)\in n_{2j+1}\times \ou}$ of finite sets of pairs from $\N\times \ou$ such that

\noindent (c) $|F_i^\al|=n_{2j_i}$ for every $(i,\al)\in n_{2j+1}\times \ou$,

\noindent (d) $(y_k^\xi)_{(k,\xi)\in F_i^\al}$ is a $(3,(n_{2j_i})^{-2})$-RIS and
$(f_k^\xi)_{(k,\xi)\in F_i^\al}$ is a separated block sequence for every $(i,\al)\in n_{2j+1}\times
\ou$.

 \noindent (e)    $|f_k^\xi(T(y_k^\xi))|, \nrm{x_k^{\xi}-y_k^\xi}\le n_{2j_i}^{-3}$ for  every
$(k,\xi)\in F_i^\al $ and every $(i,\al)\in n_{2j+1}\times \ou$.

\noindent (f) $2j_i=\sig_\varrho(\phi_0^\al,m_{2j_0},p_0,\dots, \phi_{i-1}^\al,
m_{2j_{i-1}},p_{i-1})$ for every $(i,\al)\in n_{2j+1}\times \ou$ where   $p_i$ is an integer such
that
\begin{align*} p_i\ge & \max\{p_0,\dots,p_{i-1},n_{2j+1}^2,p_\varrho(\bigcup_{k<i}(\supp
\phi_k^\al\cup\supp t_k^\al)),|\supp t_{i-1}^\al| n_{2j+1}^2\},
\end{align*}
and where
\begin{align*}
\phi_i^\al=&\frac1{m_{2j_i}}\sum_{(k,\xi)\in F_i^\al}f_k^\xi,\\
t_i^\al=& \frac{n_{2j_i}}{m_{2j_i}}\sum_{(k,\xi)\in F_i^\al}y_k^\xi.
\end{align*}
and $\sig_\varrho$ and $p_\varrho$ are the coding and the $\varrho$-number, respectively,
introduced in \cite{arg-lop-tod}. Set also
\begin{align*}
z_i^\al=& \frac{n_{2j_i}}{m_{2j_i}}\sum_{(k,\xi)\in F_i^\al}x_k^\xi.
\end{align*}
 Now again using Corollary \ref{oeuwtu3432423fgn} we can find countable ordinals
$\al_0<\dots <\al_{n_{2j+1}-1}$ such that $(\phi_i^{\al_i})_{i<n_{2j+1}}$ is a separated block sequence. It
follows that the sequence $((t_i^{\al_i},\phi_i^{\al_i}))_{i<n_{2j+1}}$ is a $(n_{2j+1}^{-2},j)$-dependent
sequence, a slightly variation of the notion of $(0,j)$-dependent sequence used in \cite{arg-lop-tod}
(Definition 5.22; see also the proof of Proposition 5.24). The only change is that now one has that
$|\phi_i^{\al_i}(t_i^{\al_i})|\le 1/n_{2j+1}^2$ instead of zero. It follows that
\begin{equation}
\nrm{\frac{1}{n_{2j+1}}\sum_{i<n_{2j+1}} t_i^{\al_i} }\le \frac{1}{m_{2j+1}^2}.
\end{equation}
Hence, setting $z=(1/n_{2j+1})\sum_{i<n_{2j+1}} z_i^{\al_i}$, one has that
\begin{equation}
\label{joertjitjioe}\nrm{z}\le \frac{2}{m_{2j+1}^2}.
\end{equation}
By the other  hand, since $\phi=(1/m_{2j+1})\sum_{i<n_{2j+1}}\phi_i^{\al_i}$ is in $K$, it follows
that
\begin{equation}
\label{joertjitjioe1}\nrm{T(z) }\ge \phi(T(z))\ge \frac{\vep}{m_{2j+1}}.
\end{equation}
Putting \eqref{joertjitjioe} and \eqref{joertjitjioe1} together one gets
$$\frac{\vep}{m_{2j+1}}\le \nrm{T(z)}\le \nrm{T}\nrm{z}\le \frac{2\nrm{T}}{m_{2j+1}^2},$$
and this is contradictory with the choice of $j$.
\fprue
Now for each countable ordinal $\al$, let $\la_k^{\al}=\la_k^\al(T,\bsy x,\bsy y)\in \R$ be such that
$$d(T(x_{k}^\al),\R x_k^\al)=\nrm{T(x_k^\al)-\la_k^\al x_k^\al}.$$
\lema\label{kjnfjkggf}
For all but countably many $\al<\ou$, the numerical sequence $(\la_k^\al)_k$ is convergent.
\flema
\prue
Otherwise, using Lemma \ref{ojrtreuret}, one can find two real numbers $\de<\vep$, an uncountable
set $A\con \ou$, for each $\al\in A$ two infinite disjoint subsets $L_\al$ and $R_\al$ of $\N$, and
a block sequence $(f_k^\al)_{(k,\al)\in (L_\al\cup R_\al)\times A}$ in $K$ such that

\noindent (i) $\ran f_k^\al\con \ran y_k^\al$, $f_k^\al(y_k^\al)=1$ for every $(k,\al)\in
(L_\al\cup R_\al)\times A$.

\noindent (ii) For every $\al\in A$ one has that   $f_k^\al(T(x_k^\al))<\de$ if $k\in L_\al$, and
$f_k^\al(T(x_k^\al))>\vep$ if $k\in R_\al$.

Let $j\in \N$ be such that $(\vep-\de) m_{2j+1}\ge 4 \nrm{T}$.  We find, as in the proof of the
previous Lemma \ref{ojrtreuret} a sequence  $(F_i^\al)_{(i,\al)\in n_{2j+1}\times \ou}$ of finite
sets of pairs from $\N\times \ou$ such that  (c), (d),  and (f) as there holds, and also

\noindent (e')  For  every $(k,\xi)\in F_i^\al $ one has that $\nrm{x_k^{\xi}-y_k^\xi}\le
n_{2j_i}^{-3}$  for  every $(i,\al)\in n_{2j+1}\times \ou$, and    $f_k^\xi(T(x_k^\xi))\ge \vep$ if
$i$ odd  and $f_k^\xi(T(x_k^\xi))\le \de$ if $i$ even.

We set also $\phi_i^\al$, $u_i^\al$ and $z_i^\al$ as there.  Let $\al_0<\dots <\al_{n_{2j+1}-1}$ be
such that $(\phi_i^{\al_i})_{i<n_{2j+1}}$ is a separated block sequence. It follows that the
sequence $((t_i^{\al_i},\phi_i^{\al_i}))_{i<n_{2j+1}}$ is a $(1,j)$-dependent sequence (see
Definition 3.3 in \cite{arg-lop-tod}), hence
\begin{equation}
 \nrm{\frac{1}{n_{2j+1 }}\sum_{i<n_{2j+1}}(-1)^i u_i^{\al_i}}\le \frac1{m_{2j+1}^{2}},
\end{equation}
and so by the property (e') one has that
\begin{equation}
\label{oiejhotiuragj} \nrm{z}\le \frac2{m_{2j+1}^{2}},
\end{equation}
where $z=(1/n_{2j+1 })\sum_{i<n_{2j+1}}(-1)^i z_i^{\al_i}$. One also has, setting
$\phi=(1/m_{2j+1})\sum_{i<n_{2j+1}}\phi_i^{\al_i}$, that
\begin{equation}\label{oiejhotiuragj1}
\nrm{T(z)}\ge
|\phi(Tz)|=|\frac{1}{m_{2j+1}n_{2j+1}}\sum_{i<n_{2j+1}}(-1)^i\phi_i^{\al_i}(T(z_i^{\al_i})) |\ge
\frac{\vep-\de}{2m_{2j+1}}.
\end{equation}
From \eqref{oiejhotiuragj} and \eqref{oiejhotiuragj1} one easily gets a contradiction with the
choice of $j$.
\fprue
For every $\al<\ou$  let (if exists) $\la_\al=\la_\al(T,\bsy x,\bsy y)=\lim_{k\to \infty }
\la_k^\al$.
\cor\label{ohrhtughr}
There is a real number $\la=\la(T,\boldsymbol{x},\boldsymbol{y})$ such that $\la_\al=\la$ for all
but countably many $\al$.
\fcor

\prue
Otherwise, there are two reals $\de<\vep$ such that both sets $A_0=\conj{\al<\ou}{\la_\al<\de}$ and
$A_1=\conj{\al<\ou}{\la_\al>\vep}$ are uncountable. Find for every $(k,\al)\in \N\times \ou$ a countable
ordinals $\be(k,\al)$ ($(k,\al)\in \N\times \ou$) such that

\noindent (i) $\be(k,\al)\in A_i$  if $k$ is equal to $i$ $\mod 2$.

\noindent (ii) $((x_k^{\be(k,\al)},y_k^{\be(k,\al)}))_{(k,\al)\in \N\times \ou}$ is a deciding pair
$(\bsy z,\bsy t)$ for $T$.

It follows that $(\la_k^\al(\bsy z,\bsy t))_k$ is never converging, contradicting Lemma
\ref{kjnfjkggf}.
\fprue
\cor
The scalar $\la(T,\boldsymbol{x},\boldsymbol{y})$ is independent of $\boldsymbol{x}$ and
$\boldsymbol{y}$. We call it  $\la(T)$.
\fcor
\prue
Fix two deciding pairs $(\bsy x,\bsy y)$ and $(\bsy z,\bsy u)$ for $T$. It is easy to define a
third one $(\bsy v,\bsy w)$ for $T$ such that  the sets
\begin{align*}
\{\al<\ou\,& :\,(\exists \be<\ou) (\forall k\in \N) \, (x_k^\al,y_k^\al)=(z_k^\al,u_k^\al)\} \\
\{\al<\ou\, & : \,(\exists \be<\ou) (\forall k\in \N) \, (v_k^\al,w_k^\al)=(z_k^\al,u_k^\al)\}
\end{align*}
are
both uncountable. It follows that $\la(T,\bsy x,\bsy y)=\la(T,\bsy v, \bsy w)=\la(T,\bsy z,\bsy
u)$.
\fprue

\teor\label{ouetiuoeugov}
The mapping $\la:\mathcal{L}(X,\eqs)\to \R$ which sends $T$ to $\la(T)$ is  a bounded linear
functional  whose kernel $\ker(\la)$   is equal to the family of all $\ou$-singular operators from
$X$ into $\eqs$.
\fteor

\prue It is obvious from the definition of $\la(T)$ that $|\la(T)|\le \nrm{T}$. We now show the linearity of $\la$.
 It is easy to see that $\la( \mu T)=\mu \la (T)$. Let
us prove now that $\la(T_0+T_1)=\la(T_0)+\la(T_1)$.
  Let $(\boldsymbol{x},\boldsymbol{y})$ be a deciding pair for both $T_0$ and $T_1$ (See Remark \ref{ojgtiorufjnw}). Let
$\la_k^{\al,0},\la_k^{\al,1},\la_k^\al\in \R$ be such that $d((T_0+T_1)(x_k^\al),\R
x_k^\al)=\nrm{(T_0+T_1)(x_k^\al)-\la_k^\al x_k^\al}$, and $d(T_i(x_k^\al),\R x_k^{\al})=\nrm{T
(x_k^{\al})-\la_k^{\al,i} x_k^{\al}}$, for $i=0,1$. It follows from Lemma \ref{ojrtreuret}, applied
to $T_0,T_1$ and $T_0+T_1$, that
\begin{equation}
\label{etjirutirgffg}\lim_{k\to \infty}(\la_k^{\al,0}+\la_k^{\al,1})-\la_k^{\al}=0,
\end{equation}
for all but countably many $\al$. The desired result now follows    from \eqref{etjirutirgffg} and
Lemma \ref{kjnfjkggf}.

Now we prove that $\ker(\la)$ is the family of the $\ou$-singular operators. Suppose first that
$\la(T)=0$.   We are going to show that $T$ is not isomorphism when restricted to a non-separable
subspace of $X$. To do this, let $\vep>0$, and let $Z$ be a non-separable subspace of $X$. Let
$(\boldsymbol{x},\boldsymbol{y})$ be any deciding pair for $T$ with $\boldsymbol{x}\con Z$. Since
$\la(T)=0$, by previous Lemma \ref{ojrtreuret}  we can find $(k,\al)\in  \N\times\ou$ such that
$\nrm{T(x_k^{\al})}<\vep$, as desired.

Finally, suppose that $T$ is a $\ou$-singular operator. Our intention is to provide a deciding pair
$(\bsy x, \bsy y)$ such that $\nrm{T(x_k^\al)}\le 2^{-k}$ for every $(k,\al)$. Then $\la_\al(T,\bsy
x,\bsy y)=0$ for every $\al$, and so $\la(T)=0$.
\clam
For every non-separable  $X_0\hookrightarrow X$ every $\vep>0$ and every $k$ there  are two
normalized vectors $x$ and $y$ such that  $x\in X_0$, $y$ is a  $2-\ell_1^k$-average,
$\nrm{T(x)}\le \vep$ and $\nrm{x-y}\le \vep$.
\fclam
It is easy to find the desired deciding pair $(\bsy x,\bsy y)$ from a simple use of the previous
claim. Now we pass to give a proof of the claim. Fix all data. Let $l\in \N$ be  the result of the
application of Proposition \ref{dfmlejtilejiew} to our fixed $k$, and let $\de=\vep/2l$. Since $T$
is $\ou$-singular, one can find two normalized uncountable sequences $(x_\al)_{\al<\ou}$ and
$(y_\al)_{\al<\ou}$ such that

\noindent (a) $x_\al\in X_0$ for every $\al<\ou$.

\noindent (b) $\bsy  y=(y_\al)_{\al<\ou}$ is a block sequence, and there an uncountable block
sequence $(f_\al)_{\al<\ou}$ in $K$ biorthogonal to $\bsy y$.

\noindent (c) $\nrm{x_\al-y_\al}, \nrm{T(x_\al)}\le \de$.

It follows from Corollary \ref{oeuwtu3432423fgn} when applied to $(\supp f_\al)_{\al<\ou}$ that
there is $F\con \ou$ of size $l$ such that $(f_\al)_{\al\in F}$ is a separated sequence. So, by
Proposition \ref{dfmlejtilejiew} one can find a $2-\ell_1^k$ in $\langle y_\al \rangle_{\al\in F}$
and then, using property (c) above, the counterpart $x\in \langle x_\al \rangle_{\al\in F}$ so that
$x$ and $y$ fulfills the desired conditions.
\fprue

From Theorem \ref{ouetiuoeugov} one easily gets the main conclusion of this section which gives us
the description of the spaces of operators $\mc L(X,\eqs)$ where $X$ is an arbitrary closed
infinite dimensional subspace of $\eqs$.
\cor
Every bounded operator $T:X\to \eqs$ from a closed subspace of $\eqs$ into $\eqs$ can be expressed
as the sum
$$T=\la(T)i_{X,\eqs}+S$$
where  $S$  is a $\ou$-singular operator.\qed
\fcor

\cor
The space $\eqs$ is $\ou$-hereditarily indecomposable, and therefore it contains no uncountable
unconditional basic sequences.
\fcor
\prue
This follows from Proposition \ref{jyhruthr}.
\fprue

It is well known  that the class of strictly singular operators on a Banach space $X$ is a closed
ideal of the Banach algebra $\mathcal{L}(X)$ of all bounded operators from $X$ into $X$. We do not
know if the same is true for the class of $\ou$-singular operators. However, we have the following.
\cor
If $X$ is a closed infinite dimensional subspace of $\eqs$, then   the family $\mc S_{\ou}(X)$ of
$\ou$-singular operators in $ \mathcal{L}(X)$ forms a closed ideal in the Banach algebra
$\mathcal{L}(X)$.
\fcor
\prue
We show that $\la:\mc{L}(X)\to \R $, formally defined by $\la \mapsto \la(T)=\la(i_{X,\eqs}\circ
T)$,  is a bounded operator between Banach \emph{algebras}. From this one easily gets the desired
result. Observe that we did almost all the work in Theorem \ref{ouetiuoeugov}.   It remains to show
that $\la(T_1\circ T_0)=\la(T_1) \la(T_0)$.   Fix a deciding pair $(\bsy x,\bsy y)$ for both $T_0$
and $T_1$. It follows that for all but countably many $\al<\ou$ one has that
\begin{align*}
\lim_{k\to \infty}\nrm{(T_1\circ T_0)(x_k^{\al})- \la(T_1\circ T_0)x_k^\al}= & 0 \\
\lim_{k\to \infty}\nrm{(T_1\circ T_0)(x_k^{\al})-T_1(\la(T_0)x_k^\al)}= & 0 \\
\lim_{k\to \infty}\nrm{T_1(\la(T_0) x_k^{\al})-\la(T_1)\la(T_0)x_k^\al}= & 0,
\end{align*}
and so $\la(T_1\circ T_0)=\la(T_1)\la (T_0)$, as desired.
\fprue

In the case of $X=\eqs$ we   obtain the following slightly more informative result.
\cor\label{oieytouwyyghds}If  $T:\eqs\to \eqs$ is a bounded operator, then
$T(e_\al)=\la(T)e_\al$ for all but countably many $\al$.  It follows that $T$ is the sum of a
multiple of the identity plus an operator with separable range.
\fcor
\prue
Otherwise we can find an uncountable subset $A\con \ou$ and $\vep>0$ such that
\begin{equation} \label{ojriyojgh}\nrm{T(e_\al)-\la(T) e_\al}>\vep
\end{equation}
for all $\al\in A$ and also, as a consequence of Remark \ref{oieutiougf} (f), such that
$T(e_\al)<T(e_\be)$ for every $\al<\be$ in $A$ (here we are making an abuse of language by formally
accepting  that $0<x$ for every vector $x$). Let $\theta_A:\om\times \omega_1\to A$ be the unique
order-preserving onto mapping, and define $x_k^\al=e_{\theta_A(k,\al)}$. Then
$((x_k^\al),(x_k^\al))$ is a deciding pair for $T$. Observe that this is what it makes so peculiar
the situation $X=\eqs$.  So, by Lemma \ref{ojrtreuret}, for all but countably many $\al$ one has
that
$$\lim_{k\to \infty}\nrm{T(x_k^\al)-\la(T)x_k^\al}=0,$$ which is contradictory with
\eqref{ojriyojgh}.
\fprue
The following result shows that, among the bounded operators from $\eqs$ into $\eqs$,  there is no
distinction between the notions of $\ou$-singular operators and the notion of operators with
separable range.
\cor
The following are equivalent for $T\in \mathcal{L}(\eqs)$:

\noindent (a) $(T(e_\al))_{\al<\ou}$ is eventually zero.

\noindent (b) $T$ has separable range.

\noindent (c) $T$ is $\ou$-singular.
\fcor
\prue
(c) implies (a): Suppose that $T$ is $\ou$-singular. Then by Theorem \ref{ouetiuoeugov} one has
that $\la(T)=0$; so by Corollary \ref{oieytouwyyghds} it follows  that $T(e_\al)=0$ eventually.
\fprue

\end{document}